
\magnification1200

\input amstex.tex
\documentstyle{amsppt}

\hsize=12.5cm
\vsize=17cm
\hoffset=1cm
\voffset=2cm

\def\DJ{\leavevmode\setbox0=\hbox{D}\kern0pt\rlap
 {\kern.04em\raise.188\ht0\hbox{-}}D}
\def\dj{\leavevmode
 \setbox0=\hbox{d}\kern0pt\rlap{\kern.215em\raise.46\ht0\hbox{-}}d}

\def\txt#1{{\textstyle{#1}}}
\baselineskip=13pt
\def\hf{{\textstyle{1\over2}}}

\def\d{{\,\roman d}}
\def\e{\varepsilon}
\def\f{\varphi}
\def\G{\Gamma}

\def\s{\sigma}

\def\={\;=\;}

\def\zt{\zeta(\hf+it)}

\def\R{\Re{\roman e}\,} 
\def\z{\zeta}

\def\hf{{\textstyle{1\over2}}}
\def\txt#1{{\textstyle{#1}}}
\def\f{\varphi}

\font\tenmsb=msbm10
\font\sevenmsb=msbm7
\font\fivemsb=msbm5
\newfam\msbfam
\textfont\msbfam=\tenmsb
\scriptfont\msbfam=\sevenmsb
\scriptscriptfont\msbfam=\fivemsb

\font\ff=cmr8
\def\txt#1{{\textstyle{#1}}}
\baselineskip=13pt

\font\teneufm=eufm10
\font\seveneufm=eufm7
\font\fiveeufm=eufm5
\newfam\eufmfam
\textfont\eufmfam=\teneufm
\scriptfont\eufmfam=\seveneufm
\scriptscriptfont\eufmfam=\fiveeufm
\def\mathfrak#1{{\fam\eufmfam\relax#1}}

\font\tenmsb=msbm10
\font\sevenmsb=msbm7
\font\fivemsb=msbm5
\newfam\msbfam
     \textfont\msbfam=\tenmsb
      \scriptfont\msbfam=\sevenmsb
      \scriptscriptfont\msbfam=\fivemsb

  \def\rightheadline{{\hfil{\ff
  On the integral of Hardy's function}\hfil\tenrm\folio}}

  \def\leftheadline{{\tenrm\folio\hfil{\ff
   Aleksandar Ivi\'c }\hfil}}
  \def\emptyheadline{\hfil}
  \headline{\ifnum\pageno=1 \emptyheadline\else
  \ifodd\pageno \rightheadline \else \leftheadline\fi\fi}

\font\ff=cmr8
\font\teneufm=eufm10
\font\seveneufm=eufm7
\font\fiveeufm=eufm5
\newfam\eufmfam
\textfont\eufmfam=\teneufm
\scriptfont\eufmfam=\seveneufm
\scriptscriptfont\eufmfam=\fiveeufm
\def\mathfrak#1{{\fam\eufmfam\relax#1}}

\font\tenmsb=msbm10
\font\sevenmsb=msbm7
\font\fivemsb=msbm5
\newfam\msbfam
\textfont\msbfam=\tenmsb
\scriptfont\msbfam=\sevenmsb
\scriptscriptfont\msbfam=\fivemsb

 \def\e{\varepsilon}
 \def\d{\,{\roman d}}
\topmatter
\title
ON THE INTEGRAL OF HARDY'S FUNCTION
\endtitle
\author
Aleksandar Ivi\'c
\endauthor
\address
Katedra Matematike RGF-a Universiteta u Beogradu, \DJ u\v sina 7,
11000 Beograd, Serbia (Yugoslavia)
\endaddress
\keywords The Riemann zeta-function, Hardy's function, saddle point
\endkeywords
\subjclass 11 M 06
\endsubjclass
\email {\tt aivic\@matf.bg.ac.rs, ivic\@rgf.bg.ac.rs}
\endemail
\dedicatory
To the memory of A.A. Lavrik (1964-2003)
\enddedicatory
\medskip
\centerline{
{\sevenbf Archiv Mathematik 83(2004), 41-47}}
\medskip
\abstract
If $Z(t) = \chi^{-1/2}(\hf+it)\zt$ denotes Hardy's function,
where $\z(s) = \chi(s)\z(1-s)$, then it is proved that
$$
\int_0^T Z(t)\d t = O_\e(T^{1/4+\e}).
$$
\endabstract
\endtopmatter
\document

Let as usual $\z(s) = \sum_{n=1}^\infty n^{-s}\;(\s > 1)$ denote
the Riemann zeta-function, where $s = \s+it$ is a complex
variable. The aim of this note is to provide a
bound for the integral of Hardy's function
$$
Z(t) = \chi^{-1/2}(\hf+it)\zt,\quad \chi(s) = 2^s\pi^{s-1}\sin(\hf\pi s)
\G(1-s),\leqno(1)
$$
so that the functional equation for $\z(s)$ has the form
$\z(s) = \chi(s)\z(1-s)$.
Since $\chi(s)\chi(1-s) = 1$, it follows that $|Z(t)| = |\zt|$, and
that $Z(t)$ is a real-valued function of $t$. The function $Z(t)$ plays
an important r\^ole in the theory of the distribution of zeros of $\z(s)$
on the ``critical line" $\R s = \hf$ (see e.g., [1]--[3] and [5]-[6]).

\smallskip
The result on the integral of $Z(t)$ is contained in the following

\medskip
THEOREM. {\it We have}
$$
\int_0^T Z(t)\d t \;=\; O_\e(T^{{1\over4}+\e}). \leqno(2)
$$

\bigskip
{\bf Proof}. Here and later $\e$ will denote arbitrarily small, positive
constants, not necessarily the same ones at each occurrence.
To prove (2) we shall make use of the approximate functional equation
$$
Z^k(t) = 2\sum_{n\le2\tau}\rho\left({n\over\tau}\right)d_k(n)n^{-1/2}
\cos\left(t\log{\tau\over n} - {k\over2}t - {\pi k\over8}\right)
+ O(t^{{k\over4}-1}\log^{k-1}t),\leqno(3)
$$
which is valid for any fixed integer $k\ge1$ and $t\ge 2$.
In (3) we have set for brevity
$$
\tau \;=\; \left(t\over2\pi \right)^{k/2},
$$
and further notation is as follows. The function $d_k(n)$ represents
the number of ways $n$ may be represented as the product of $k$
factors ($d_1(n) \equiv 1, d_2(n) \equiv d(n)$, the number of divisors
of $n$), while $\rho(x)$ is a non-negative, smooth function
supported in $\,[0,2]\,$, such that $\rho(x) = 1$ for $0 \le x\le 1/b$
for a fixed constant $b>1$, and
$\rho(x) + \rho(1/x) = 1$ for all $x$. The author proved [4, Theorem 4.2] the
approximate functional equation for $\z^k(s)$, which gives (3) with
$x = y = \tau$, on using (1) and the asymptotic formula
$$
\chi(s) = \left({2\pi\over t}\right)^{\s+it-1/2}{\roman e}^{i(t+\pi/4)}
\cdot\left(1 + O\left(1\over t\right)\right)\qquad(0 \le \s \le 1,\;
t \ge t_0 > 0).
$$
Taking $k=1$ in (3) it follows  that
$$
\int_T^{2T}Z(t)\d t = 2\int_T^{2T}\sum_{n\le2\tau}
\rho\left({n\over\tau}\right)n^{-1/2}
\R \left\{{\roman e}^{iF(t)}\right\}\d t + O(T^{1/4}),\leqno(4)
$$
where
$$
\tau = \sqrt{t\over2\pi},\qquad F(t) = t\log{\tau\over n} - {t\over2}
- {\pi\over8}.\leqno(5)
$$
The reason that (3) was used  is that the standard approximate
functional equation for $\z(s)$ (this is the Riemann-Siegel formula,
see e.g.,  [1, Chapter 4]) has the error
term $O(t^{-1/4})$,
which is not sufficiently good to produce the bound in (2). For this
reason we resorted to (3), which is a smoothed variant of the approximate
functional equation with a sharp error term.

\smallskip
In view of (4), to prove (2) it clearly suffices to prove that
$$
I(T) :=   \int_T^{2T}\sum_{n\le2\tau}\rho\left({n\over\tau}\right)n^{-1/2}
\R \left\{{\roman e}^{iF(t)}\right\}\d t  \ll_\e T^{{1\over4}+\e}.\leqno(6)
$$
We have, in view of (5),
$$
\eqalign{I(T) &=
\sum_{n\le2\sqrt{T/\pi}}n^{-1/2}\R \left\{
\int_{T_1}^{2T}\rho\left({n\over\tau}\right){\roman e}^{iF(t)}\d t\right\}
\cr&= \sum\nolimits_1(T) + \sum\nolimits_2(T) +  \sum\nolimits_3(T) +
 \sum\nolimits_4(T) + \sum\nolimits_5(T),\cr}\leqno(7)
$$
say, where
$$
T_1 = \max\left(T,\,2\pi{\left({n\over2}\right)}^2\right),
$$
and the ranges of summation in $\sum_j(T)\;(j = 1,\ldots, 5)$
are respectively
as follows: $n\le \sqrt{T/(2\pi)} - T^{\e}$, $\sqrt{T/(2\pi)} - T^{\e}
< n \le \sqrt{T/(2\pi)} + T^{\e}$, $\sqrt{T/(2\pi)} + T^{\e} < n
\le \sqrt{T/\pi} - T^{\e}$, $\sqrt{T/\pi} - T^{\e} < n \le
\sqrt{T/\pi} + T^{\e}$ and $\sqrt{T/\pi} + T^{\e} < n \le 2\sqrt{T/\pi}$.
We have
$$
F'(t) \;=\; \log{\sqrt{t/(2\pi)}\over n},\quad F''(t) \;=\; {1\over2t}.
$$
This means that, in $\sum_1(T)$, we have
$$
F'(t) \ge \log{[\sqrt{T/(2\pi)}]\over n},
$$
hence ($m = [\sqrt{T/(2\pi)}] - n)$) by the first derivative test
(see e.g., [1, Lemma 2.1])
$$
\eqalign{
\sum\nolimits_1(T) & \ll T^{1/4}
+ \sum_{\hf\sqrt{T/(2\pi)}< n \le \sqrt{T/(2\pi)} - T^{\e}}
{1\over\sqrt{n}\log ([\sqrt{T/(2\pi)}]/n)}\cr&
\ll T^{1/4}
+ T^{-1/4}\sum_{\hf\sqrt{T/(2\pi)}< n \le \sqrt{T/(2\pi)} - T^{\e}}
{n\over [\sqrt{T/(2\pi)}] - n}\cr&
\ll T^{1/4}
+ T^{1/4}\sum_{m \le \hf\sqrt{T/(2\pi)}}
{1\over m}\cr& \ll T^{1/4}\log T.\cr}
$$
An analogous bound holds also for $\sum\nolimits_5(T)$.

To evaluate the sum  $\sum_3(T)$ in (7), which contains
(for every $n$ in the range of summation) a saddle point $c$, namely the
root of $F'(c) = 0$, so that $c = c_n = 2\pi n^2$,
one may use general results in the literature which  for this purpose
(see [1], [4] and [5]). A convenient one is [5, Lemma III.2], which says that
$$
\eqalign{&
\int_a^b \f(x)\exp(2\pi if(x))\d x = {\f(c)\over\sqrt{f''(c)}}
{\roman e}^{2\pi if(c)+\pi i/4} + O(HAU^{-1})\cr&
+ O\bigl(H\min(|f'(a)|^{-1},\sqrt{A}\,\bigr) +
O\bigl(H\min(|f'(b)|^{-1},\sqrt{A}\,\bigr),\cr}
\leqno(8)
$$
if $f'(c) = 0$, $a \le c \le b$, and the following conditions hold: $f(x) \in
C^4[a,b]$, $\f(x) \in C^2[a,b]$, $f''(x) > 0$ in $[a,b]$, $f''(x) \asymp
A^{-1}$, $f^{(3)}(x) \ll A^{-1}U^{-1}$, $f^{(4)}(x) \ll A^{-1}U^{-2}$,
$\f^{(r)}(x) \ll HU^{-r}\;(r = 0,1,2)$ in $[a,b],\;0 < H,A < U, \,
0 < b - a \le U$.

We shall apply (8) with $\f(t) = \rho(n/\tau)$, $f(t) = (2\pi )^{-1}F(t)$,
$a = T_1, b = 2T, H = 1, U = T, A = T$. With $c = c_n = 2\pi n^2$ we have
that  $c \in [a,b]$ for our range of  $n$, and
furthermore $\f(c_n) = \rho(1) = 1$. Therefore the contribution of the
first term on the right-hand side of (8) will be
$$\eqalign{&
\sum_{\sqrt{T/(2\pi)} + T^{\e} < n\le \sqrt{T/\pi} - T^{\e}}
n^{-1/2}\f(c_n)(f''(c_n))^{-1/2}\R\left\{
\exp(2\pi if(c_n)+ {\txt{1\over4}}\pi i)\right\}\cr&
=  \sqrt{8}\pi
\sum_{\sqrt{T/(2\pi)} + T^{\e} < n\le \sqrt{T/\pi} - T^{\e}}
n^{-1/2}n\,\R\{\exp({\txt{1\over8}}\pi i - \pi in^2)\}\cr&
= \sqrt{8}\pi\cos\left({\pi\over8}\right)
\sum_{\sqrt{T/(2\pi)} + T^{\e} < n\le \sqrt{T/\pi} - T^{\e}}
(-1)^n n^{1/2} = O(T^{1/4}),
\cr}\leqno(9)
$$
since the last sum is, in absolute value,
$$
\le \left|\sum_{\ell\ll K}(\sqrt{K+2\ell}-\sqrt{K+2\ell-1})\right|
\ll \sum_{\ell\ll K}1/\sqrt{K} \ll \sqrt{K}\quad(K \asymp \sqrt{T}).
$$

In the $\sum_3(T)$ we have
$\sqrt{T/(2\pi)} + T^{\e} < n \le \sqrt{T/\pi} - T^{\e}$, hence
similarly to the estimation of $\sum_1(T)$, the  total
contribution of the error terms in (8) will be $\ll_\e T^{1/4+\e}$.

Finally, by using the second derivative test ([1, Lemma 2.2]), it follows
that
$$
\sum\nolimits_2(T)  + \sum\nolimits_4(T) \,\ll_\e\,
T^{\e}T^{1/4}T^{-1/4}T^{1/2}
= T^{1/2+\e}.\leqno(10)
$$
Therefore, except for the bound
in (10), we get the upper bound $O_\e(T^{1/4+\e})$ for
our integral $I(T)$ (see (6)). The reason
for the range of summation over $n$ in $\sum_3(T)$ was the structure
of the error terms in (8), namely if $a$ or $b$ is too near a
saddle point, then $\sqrt{A}$ is to be taken, which in our case is
too large to produce (2).

\medskip
To get around this obstacle, we shall employ the saddle point method
directly, taking advantage of the particular structure of the
exponential integrals in question, coupled with the summation over
$n$ in (7). The main terms will be, of course, the same ones as those
which appeared in (8), and the essential fact is the presence
of $(-1)^n$ in the summation over $n$, which accounts for massive
cancellation and leads to (2).

\medskip
Henceforth we suppose that $n$ lies in the range covered by
$\sum_j(T)\;(j=2,3,4)$ in (7), namely
$$
\sqrt{T\over2\pi} - T^\e \le n \le \sqrt{T\over\pi} + T^\e.
$$
For such $n$ let
$$
J(T,n) = [2\pi n^2 - T^{\e},\,2\pi n^2 + T^{\e}],\quad
K(T,n) = [T_1,\,2T] \,\backslash\, J(T,n).\leqno(11)
$$
In dealing with
$$
\int_{K(T,n)}\rho\left({n\over\tau}\right){\roman e}^{iF(t)}\d t
$$
we apply the first derivative test as before, obtaining after
summation over $n$ a contribution which is $\ll_\e T^{1/4+\e}$. In case
$J(T,n)$ does not entirely lie in $[T,\,2T]$, obvious modifications in the
argument are to be made. To evaluate
$$
\int_{J(T,n)}\rho\left({n\over\tau}\right){\roman e}^{iF(t)}\d t,\leqno(12)
$$
we develop first $\rho({n\over\tau})$  by Taylor's formula at
the point $2\pi n^2$. Since each derivative of
$$
\rho\left({n\over\tau}\right) = \rho\left({n\over\sqrt{t\over2\pi}}\right),
$$
as a function of $t$, decreases by a factor of $T$, and the measure
of $J(T,n)$ is $\ll_\e T^{\e}$, we first take so many terms in Taylor's
formula so that the contribution of the error term is negligible,
namely $\ll_\e T^{1/4+\e}$. The remaining integrals will be all of the
same type, with the same exponential factor, and the largest one will
be the first one, namely the one with ($c_n = 2\pi n^2$)
$$
\rho\left({n\over\sqrt{c_n\over2\pi}}\right) = \rho(1) = {1\over2},
$$
since $\rho(x) + \rho(1/x) = 1$. Then we write, by Cauchy's theorem,
$$
\int_{J(T,n)}{\roman e}^{iF(z)}\d z =
\int_{\Cal L_1}{\roman e}^{iF(z)}\d z +
\int_{\Cal L_2}{\roman e}^{iF(z)}\d z +
\int_{\Cal L_3}{\roman e}^{iF(z)}\d z,\leqno(13)
$$
say, where $\Cal L_1$ is the segment $c_n - T^{\e}
+ v{\roman e}^{-{1\over4}\pi i}$, $0 \le v \le {1\over\sqrt{2}}T^{\e}$,
$\Cal L_2$ is the segment $c_n + v{\roman e}^{{1\over4}\pi i}$,
$|v| \le {1\over\sqrt{2}}T^{\e}$, and
$\Cal L_3$ is the segment $c_n + T^{\e}
- v{\roman e}^{-{1\over4}\pi i}$, $0 \le v \le {1\over\sqrt{2}}T^{e}$.
On $\Cal L_2$ we have
$$
iF(z) = iF(c_n) + i{v^2\over2!}{\roman e}^{{1\over2}\pi i}F''(c_n)
+ i{v^3\over3!}{\roman e}^{{3\over4}\pi i}F'''(c_n)
+ i{v^4\over4!}F^{(4)}(c_n) + \cdots\,. \leqno(14)
$$
Note that
$$
v^{k}F^{(k)}(c_n) \;\ll_{k,\e}\; T^{k\e}T^{1-k} \;=\; T^{1-k+k\e}\qquad(k = 2,3,\ldots\,).
$$
Hence if we choose $K = K(\e)$ sufficiently large, then the terms
of the series in (14) for $k > K$,
on using $\exp z = 1 + O(|z|)$ for $|z| \le 1$,
will make a negligible contribution. Then we have
$$
\exp(iF(z)) = \exp(iF(c_n))\exp(-\hf v^2 F''(c_n))
\exp\left(\sum_{k=3}^K d_k v^k F^{(k)}(c_n)\right)
$$
with  $d_k = \exp((k+2){\pi i\over4})/k!$.
The last exponential factor is expanded by Taylor's series, and again
the terms of the series (with $v^k$) for $k > K$ will make a negligible
contribution. In the remaining terms we restore integration over
$v$ to the whole real line, making a very small error. Then we use
the classical integral (see e.g., the Appendix of [1])
$$
\int_{-\infty}^\infty \exp(Ax - Bx^2)\,\d x \;=\;
\sqrt{\pi\over B}\exp\left({A^2\over4B}\right)\qquad(\R B > 0).\leqno(15)
$$
By differentiating (15) as a function of $A$ we may explicitly
evaluate integrals of the type
$$
\int_{-\infty}^\infty x^{2k}\exp(-Bx^2)\,\d x\qquad(\R B > 0,\;
k = 0,1,2,\ldots\,).
$$
It transpires that the largest contribution ($= \sqrt{\pi}$) will come
from the integral with $k = 0$, which will coincide with the contribution
of the main term in (8).

It remains to deal with the integrals over $\Cal L_1$ and $\Cal L_3$
in (13), which are estimated analogously, so only the former is
considered. On $\Cal L_1$ we have
$$\eqalign{&
iF(c_n - T^{1-\e} + v{\roman e}^{-i\pi/4})\cr&
= i\Bigl\{F(c_n - T^{1-\e}) +
F'(c_n - T^{1-\e})v{\roman e}^{-i\pi/4} \cr&+
F''(c_n - T^{1-\e}){v^2\over2!}{\roman e}^{-2i\pi/4})
+  F'''(c_n - T^{1-\e}){v^3\over3!}{\roman e}^{-3i\pi/4})
+\ldots \Bigr)\Bigr\}\cr&
= F_1(v;n,T) + iF_2(v;n,T),\cr}
$$
say, with $F_1,\,F_2$ real. Then
$$
\eqalign{&
{\partial F_1(v;n,T)\over\partial v} + i{\partial F_2(v;n,T)
\over \partial v}\cr&
= {1+i\over\sqrt{2}}F'(c_n - T^{1-\e})
+ F''(c_n - T^{1-\e})v \cr&
+ {1-i\over\sqrt{2}}F'''(c_n - T^{1-\e}){v^2\over2!} + \ldots\,.\cr}
$$
Therefore we find that
$$
{\partial F_2(v;n,T)\over \partial v} \;\gg_\e\; T^{-\e},
$$
hence by the first derivative test the total contribution of the integral
over $\Cal L_1$ is seen to be $\ll_\e T^{1/4+\e}$.
This finishes the proof of (2). However, the true order of the integral
of $Z(t)$ remains elusive. In particular, it would be of interest to
find an omega result for this quantity. Is it true that perhaps
$$
\int_0^T Z(t)\d t \;=\; \Omega(T^{1/4})\quad(= \Omega_\pm(T^{1/4}))?
\leqno(16)
$$
If yes, then the result of the Theorem would be (up to the factor ``$\e$")
best possible. The reason that (16) seems plausible is that
$T^{1/4}$ is the order  of the terms coming from
the saddle points (see (8)), and in the evaluation of
exponential integrals  one usually expects the saddle points
to produce the largest contribution.

\bigskip {\bf Note}. The revisions to the published version, made in July 2009,
involve the corrections of some misprints. In the meantime 
M.A. Korolev, ``{\it On the primitive of the Hardy function $Z(t)$}'',
Dokl. Math. 75, No. 2, 295-298 (2007); translation from Dokl. Akad. Nauk,
Ross. Akad. Nauk 413, No. 5, 599--602 (2007), proved (16). Another proof is
to be found in a forthcoming work of M. Jutila.
\vfill
\eject\topskip2cm

\bigskip\bigskip
\Refs
\bigskip

\item{[1]} A. Ivi\'c, The Riemann zeta-function, John Wiley \&
Sons, New York, 1985.

\item{[2]} A. Ivi\'c, On a problem connected with zeros of $\zeta(s)$
on the critical line, Monatshefte Math.  {\bf104} (1987), 17-27.

\item{[3]} A. Ivi\'c and M. Jutila, Gaps between consecutive zeros of the
    Riemann zeta-function, Monatshefte Math. {\bf105} (1988), 59-73.

\item{[4]} A. Ivi\'c, The mean values of the Riemann zeta-function, Tata
	Institute of Fundamental Research, Lecture Notes {\bf82},
    Bombay 1991 (distr. Springer Verlag, Berlin etc.).

\item{[5]} A.A. Karatsuba and S.M. Voronin, The Riemann zeta-function,
Walter de Gruyter, Berlin etc., 1992.

\item{[6]} A.A. Lavrik,  Uniform approximations and zeros in short
intervals of the derivatives of the Hardy function, Soviet Math. Dokl.
{\bf40} (1990), 20-22.

\vskip3cm
\endRefs

\enddocument

\bye